\documentclass[12pt,leqno]{amsart}
\usepackage{amssymb,amsmath}
\begin{document}
\renewcommand{\Large}{\large}
\renewcommand{\LARGE}{\large\bf}

\newtheorem{theorem}{Theorem}
\newtheorem{proposition}[theorem]{Proposition}
\newtheorem{fact}[equation]{}
\newtheorem{corollary}[theorem]{Corollary}
\newtheorem{definition}[theorem]{Definition}
\newtheorem{remark}[theorem]{Remark}
\newtheorem{lemma}[theorem]{Lemma}
\newtheorem{maintheorem}{Main Theorem}

\newenvironment{pf}{\noindent {\it Proof:}\/ }
{\par \par \medskip\par }

\renewcommand{\Bbb}[1]{{\bf #1}}
\newcommand{\supertext}[1]
        {{^{\text{\scriptsize #1}}}}
\newcommand{\set}[1]{\{\,#1\,\}}
\newcommand{\pair}[2]{\langle #1,#2 \rangle}
\newcommand{\matr}[4]{\left[
\begin{array}{rr} #1 & #2 \\ #3 & #4
                    \end{array} \right]}
\newcommand{\smatr}[4]{\left[
\begin{array}{rr}
{\scriptstyle #1} & {\scriptstyle #2}\\
{\scriptstyle #3} & {\scriptstyle #4}
\end{array} \right]}
\newcommand{\cal}{\mathcal}
\newcommand{\ov}{\overline}
\newcommand{\Disc}{\operatorname{Disc}}
\newcommand{\bbQ}{{\Bbb Q}}
\newcommand{\bbC}{{\Bbb C}}
\newcommand{\cM}{{\cal M}}
\newcommand{\bbR}{{\Bbb R}}
\newcommand{\PGL}{\operatorname{PGL}}
\newcommand{\Qp}{{\Bbb Q}_{p}}
\newcommand{\Vd}{V^{(d\mbox{\hspace*{0.08em}})}}
\newcommand{\VBd}{V_{B}^{(d\mbox{\hspace*{0.08em}})}}
\newcommand{\VBp}{V_{B}^{(p\mbox{\hspace*{0.08em}})}}
\newcommand{\VBq}{V_{B}^{(q\mbox{\hspace*{0.08em}})}}
\newcommand{\tVBp}{\tilde{V}_{B}^{(p\mbox{\hspace*{0.08em}})}}
\newcommand{\Vp}{V^{(p\mbox{\hspace*{0.08em}})}}
\newcommand{\Gq}{G^{(q\mbox{\hspace*{0.08em}})}}
\newcommand{\Gqe}{G^{(q\mbox{\hspace*{0.08em}}),e}}
\newcommand{\Vq}{V^{(q\mbox{\hspace*{0.08em}})}}
\newcommand{\tVp}{\tilde{V}^{(p\mbox{\hspace*{0.08em}})}}
\newcommand{\tVB}{\tilde{V}_{B}}
\newcommand{\Md}{M^{(d\mbox{\hspace*{0.08em}})}}
\newcommand{\Bpint}{B^{\,p{\rm\mbox{-}int}}}
\newcommand{\Bqint}{B^{\,q{\rm\mbox{-}int}}}
\newcommand{\Bpintx}{B^{\,p{\rm\mbox{-}int},\times}}
\newcommand{\bbZ}{{\Bbb Z}}
\newcommand{\bbZp}{{\bbZ}_{p}}
\newcommand{\bbZup}{\bbZ^{{\rm (}p{\rm )}}}
\newcommand{\bbZupx}{{\bbZ}^{{\rm (}p{\rm )}\times}}
\newcommand{\cD}{{\cal D}}
\newcommand{\cDup}{\cD^{{\rm (}p{\rm )}}}
\newcommand{\cDupx}{\cD^{{\rm (}p{\rm )}\times}}
\newcommand{\N}{\operatorname{N}}
\newcommand{\valp}{\operatorname{val}_{p}}
\newcommand{\Gammaplus}{\Gamma_{+}}
\newcommand{\tGammaplus}{\tilde{\Gamma}_{+}}
\newcommand{\Gammaplusq}{\Gamma_{+,q}}
\newcommand{\Gammaplusp}{\Gamma_{+,p}}
\newcommand{\tGammaplusq}{{\tilde{\Gamma}}_{+,q}}
\newcommand{\G}{\operatorname{G}}
\newcommand{\Mdp}{\Md/\bbZp}
\newcommand{\tGammaplusp}{{\tilde{\Gamma}}_{+,p}}
\newcommand{\Mat}{\operatorname{Mat}}
\newcommand{\F}{\operatorname{F}}
\newcommand{\cOK}{{\cal O}_{K}}
\newcommand{\GL}{\operatorname{GL}}
\newcommand{\cP}{{\cal P}}
\newcommand{\cPunr}{\cP^{\rm unr}}
\newcommand{\tGamma}{\tilde{\Gamma}}
\newcommand{\SL}{\operatorname{SL}}
\newcommand{\bbFp}{{\Bbb F}_{p}}
\newcommand{\e}{\emptyset}
\newcommand{\Pic}{\operatorname{Pic}}
\newcommand{\Picr}{\Pic ^{r}}
\newcommand{\Qq}{{\Bbb Q}_{q}}
\newcommand{\Ql}{{\Bbb Q}_{\ell}}
\newcommand{\Ftwo}{{\Bbb F}_{2}}
\newcommand{\Ffour}{{\Bbb F}_{4}}
\newcommand{\Pone}{\mathbb{P}^{1}}
\newcommand{\VM}{V_{\cM}}
\newcommand{\M}{\operatorname{M}}
\newcommand{\MGplus}{{\rm M}_{\Gammaplus}}
\newcommand{\Zp}{{\Bbb Z}_{p}}
\newcommand{\Fp}{{\Bbb F}_{p}}
\newcommand{\Fpr}{{\Bbb F}_{p^{r}}}
\newcommand{\W}{\operatorname{W}}
\newcommand{\Norm}{\operatorname{Norm}}
\newcommand{\bk}{{\bf k}}

\title[Atkin-Lehner Quotients of Shimura Curves]
{On Atkin-Lehner Quotients
of Shimura Curves}

\author{Bruce W. Jordan}
\address{Department of Mathematics, Box G-0930\\
Baruch College, CUNY\\
17 Lexington Avenue\\
New York, NY  10010  USA}
\thanks{The first author was partially supported by grants from the NSF
and PSC-CUNY}
\author{Ron Livn\'{e}}
\address{Mathematics Institute\\
The Hebrew University of Jerusalem\\
Givat Ram, Jerusalem  91904  Israel}
\email{rlivne@sunset.ma.huji.ac.il}
\thanks{Both authors were supported by a joint Binational Israel-USA
Foundation grant}
\date{February 4, 1998}
\keywords{Shimura curves, Atkin-Lehner involutions,
Shafarevich-Tate groups, $p$\/-adic uniformization}
\subjclass{11G18,11G20,14G20,14G35,14H40}
\begin{abstract}
We study the \v{C}erednik-Drinfeld $p$\/-adic uniformization
of certain Atkin-Lehner quotients of Shimura curves over $\bbQ$.
We use it to determine over which local fields they have rational
points and divisors of a given degree.
Using a criterion of Poonen and Stoll we show that the 
Shafarevich-Tate group of their jacobians is not
of square order for infinitely many cases.
\end{abstract}
\maketitle

In  \cite{PSt} Poonen and Stoll have shown that if the
Shafarevich-Tate group of a principally polarized
abelian variety $A$ over a global field $K$ is finite,
then its order can be twice a square as well as a square.
They call $A/K$ even if the quotient of its Shafarevich-Tate
group by the maximal divisible subgroup has order a square and 
odd otherwise.  They prove
\cite[Corollary 10]{PSt} that the jacobian of a curve
$X/K$ of genus $g$ is odd if and only if the number of places
$v$ of $K$ where $X$ fails to have a $K_v$-rational divisor
of degree $g-1$ is odd.  They also show that infinitely
many hyperelliptic jacobians over $\bbQ$ are odd for every
even genus.

Let $B$ be an indefinite rational quaternion division
algebra of discriminant $\Disc B$.  Let $V_{B}/\bbQ$
be the Shimura curve corresponding to a maximal order
in $B$. In \cite{JL1} and \cite{JL2} we determined
when Shimura curves $V_B$ have rational points, rational divisors
of a given degree, and rational divisor classes of a given degree
over any local field.  The key ingredient of our analysis was
the explicit $p$-adic uniformization of these curves
given by Drinfeld \cite{Dr}. The criterion of Poonen and Stoll
then immediately yields that the jacobian of $V_B/\bbQ$ is even
(\cite[Theorem 23]{PSt}).  Atkin-Lehner involutions
$w_{d}$ for $d|\Disc B$ act on the Shimura curves $V_B$.
Set $\VBd = \VBd/\bbQ=V_{B}/w_{d}$. We will show that these
Atkin-Lehner quotients $\VBd/\bbQ$, in contrast to $V_{B}/\bbQ$,
yield examples of odd nonhyperelliptic jacobians
of arbitrarily high genus. For simplicity we consider quaternion
algebras $B$ in which exactly two odd primes ramify and study
the Atkin-Lehner quotient $\VBp/\bbQ$ where $p$ is a prime. For
a related study of more general Shimura curves over totally
real fields see \cite{JLV}. In this paper
we will prove the following:
\begin{theorem}
\label{maintheorem}
Let $p\neq q$ be
primes with $p\equiv 5\!\mod 24$ and
$q\equiv 5\!\mod 12$.  Furthermore suppose that $p$
is not a square modulo $q$.  Let $B$ be the
indefinite rational quaternion algebra with $\Disc B =pq$.
Then the jacobian of $\VBp/\bbQ$ is  odd.
\end{theorem}
The genus of the curves $\VBp$ in Theorem \ref{maintheorem}
 goes to infinity with
$\Disc B$ and we show at the end of Section 2 that
only finitely many of these curves can be hyperelliptic.
\section{Local diophantine properties of Atkin-Lehner Quotients}

Let $B$ be an indefinite rational quaternion division algebra
of discriminant $pq$, where $p$ and $q$ are distinct odd primes.
Put $\Vp/\bbQ = V_{B}/w_{p}$ and $\Vq/\bbQ = V_{B}/w_{q}$.
 For a curve $X/\bbQ$, an
extension field $K/\bbQ$, and an integer $r$ we
denote by $\Picr X(K)^{+}$ the set of divisor classes on $X$
of degree $r$ containing a $K$-rational divisor.
Let $\Bpint /\bbQ$ be the definite  quaternion algebra obtained
from $B$ by interchanging local invariants at $p$ and $\infty$;
the reduced discriminant of $\Bpint$ is $q$.  For
$\alpha,\,\beta\in\bbQ$ denote by $B(\alpha, \beta )$ the rational
quaternion algebra with basis $1$, $i$, $j$, and $k$ satisfying
$i^2=\alpha$, $j^2=\beta$, and $ij=-ji=k$.  It is ramified at a prime
$p\leq\infty $ if and only if the
Hilbert symbol $(\alpha,\beta)_{p}=-1$.
\begin{theorem}
\label{local}
{\rm 1.}\hspace*{1em}$\Pic^{1}\Vp (\bbR)^{+}\neq \e$ if and only
if $\bbQ (\sqrt{p})$ splits $B$.\\
{\rm 2.}\hspace*{1em}$\Pic^{1}\Vp (\Qp)^{+}\neq \e$.\\
{\rm 3.}\hspace*{1em}$\Pic^{1}\Vq (\Qp)^{+}\neq \e$ if and only if
\begin{equation}
\label{criterion}
\Bpint\cong B(-1,-pq )\qquad\mbox{or}\qquad \Bpint\cong B(-p,-q)\ .
\end{equation}
\end{theorem}
\begin{pf}
By applying Shimura's explicit law for the action of complex conjugation
on $V_{B}(\bbC )$, Ogg proved in \cite{Ogg} that $\Vp (\bbR )\neq \e$
if and only if $\bbQ (\sqrt{p})$ splits $B$.  It is obvious that
$\Vp (\bbR )\neq \e$ if and only if $\Pic ^{1}\Vp (\bbR )^{+}\neq \e$,
proving (1).

For (2) fix a maximal
$\bbZ [1/p]$-order $\cD\subseteq \Bpint$, so that
$\bbZ [1/p]^{\times}\subseteq \cD^{\times}$.
Let $\N$ be the
reduced norm of $\Bpint$ and let $\valp$ be the $p$-adic valuation.
Define
$$
\tGammaplus  = \{\gamma \in \cD^{\times}|\valp(\N(\gamma))\mbox{ is even}\}
\qquad\mbox{\rm and}\qquad
\Gammaplus  = \tGammaplus/ \bbZ [1/p]^{\times} \ .
$$
Fix an identification of $\Bpint\otimes\Qp$ with $\Mat_{2\times 2}(\Qp )$.
Then $\tGammaplus$ is a subgroup of
$\GL_{2}(\Qp )$ containing $p=p{\rm Id}_{2\times 2}$ and
$\Gammaplus$ is a discrete cocompact
subgroup of $\PGL_{2}(\Qp )$.

Let $\MGplus/\Zp$ denote the Mumford curve uniformized by
$\Gammaplus\subseteq
\PGL_{2}(\Qp )$.  Drinfeld proves that the curve
$V_{B}\times_{\bbQ}\Qp$ is the generic fiber of the $w_{p}$-Frobenius
quadratic twist of $\MGplus/\Zp$, see \cite[Sect. 4]{JL1}.  Hence
the quotient $\Vp$ is Mumford uniformized (without a twist)
and it has a model over $\Zp$ whose reduction consists of
$\Fp$-rational $\Pone$'s crossing at $\Fp$-rational points.
In particular over any extension $\Fpr$, $r\geq 2$, there are
smooth $\Fpr$-rational points.
By Hensel's lemma (see, e.g., \cite[Lemma 1.1]{JL1}) there
are points rational over the Witt vectors $\W (\Fpr )$, $r\geq 2$.
The trace to $\Qp$ of such a point is in $\Picr \Vp (\Qp )^{+}$.
A linear combination using relatively prime $r$'s (say $r=2,\,3$)
gives then a $\Qp$-rational divisor of degree 1, proving (2).

For (3), we have $\Pic^{2}V_{B}(\Qp )^{+}\neq \e$
by \cite[Theorem 2a]{JL2}, so also
$\Pic^{2}\Vq (\Qp )^{+}\neq \e$.  Hence
$\Pic^{1}\Vq (\Qp )^{+}\neq \e$ if and only if there is a $\Qp$-rational
effective divisor of odd degree on $\Vq$.  Such a divisor is the sum
of traces to $\Qp$ of points in extension fields of $\Qp$ at least
one of which must be of odd degree over $\Qp$.  Hence we need to
show that the condition (\ref{criterion}) holds if and only if
$\Vq (K)\neq \e$ for some extension $K/\Qp$ of odd degree.
(In fact our proof will show that condition (\ref{criterion}) is equivalent
to $\Vq (\Qp )\neq \e$ .)

Let $\Delta$ be the tree of $\SL (2, \Qp )$;
it is a homogeneous tree of regularity $p+1$.  The dual graph
(see \cite{Kur}, \cite {JL1}) of the
special fiber of $\MGplus$ is given by $G=\Gammaplus
\backslash \Delta$. Let $\gamma_{q}$ be an element of the normalizer
 of $\tGammaplus$ in $\cD^{\times}$ of
reduced norm $q$.
Set
$$
\tGammaplusq=\langle \tGammaplus , \gamma_{q}\rangle\subseteq \GL(2,\Qp )\qquad
\mbox{\rm and}\qquad
\Gammaplusq = \tGammaplusq/\bbZ[1/p]^{\times}\ .
$$
From \cite{Dr} it follows that $\gamma_{q}$ induces the action of
$w_{q}$ on $\MGplus$.
Taking the quotients of $\MGplus$ and $V_{B}$ by $w_{q}$ we see
that $\Vq\times_{\bbQ}\Qp$
is the $w_{p}$-Frobenius quadratic twist
of the Mumford curve ${\rm M}_{\Gammaplusq}$ uniformized by
$\Gammaplusq\subseteq \PGL (2,\Qp )$.  This gives a model of
$\Vq\times_{\bbQ}\Qp$ over $\Zp$ with dual graph
$\Gq=w_{q}\backslash G$ and Frobenius acting as $w_{p}$.
Let $K/\Qp$ be an extension with ring of integers $\cOK$ and residue field
$\bk$. Put
$e=e(K/\Qp )$ and $f=f(K/\Qp )$.
By base change, $\Vq\times_{\Qp}K$ has an $\cOK$-model having dual
graph $\Gq$ and Frobenius action $w_{p}^f$. Moreover, the lengths
of edges are multiplied by $e$, see \cite[Proposition 3.4]{JL1}.

The graph $G$ is bipartite with  $w_{p}$ interchanging
the even and the odd vertices; see, e.g., the discussion preceding
Proposition 4-4 in \cite{Kur}.  This property is inherited by the
quotient graph $\Gq$.
We now suppose that $K/\Qp$ is an extension of odd degree.
Then Frobenius acts on $\Gq$ by
$w_{p}$, which can fix no vertex as it interchanges the even and odd 
vertices of the bipartition. 
Hence no component of the
special fiber is rational over $\bk$, and a point $x$ of
$\Vq (K)$ must specialize to a  singularity of the special fiber.
The length $\ell (r)$ of the edge $r$ of $\Gq$ corresponding to
$x$ determines the nature of the singularity at $x$, and is
the order of the stabilizer in $\Gammaplusq$ of any preimage
$\tilde{s}$ of it in $\Delta$.
As in \cite[Theorem 5.2]{JL1}, $\Vq (K)\neq \e$
if and only if there is an oriented edge $r$ of
$\Gq$ such that $\ell (r)$ is even and $w_{p}(r)$ is the opposite
edge  $\overline{r}$.
Let $s$ be an edge of $G$ above such an $r$.  Then we see that\\

%\begin{center}
\begin{tabular}{l}
1)\hspace*{1em}Either $\ell (s)$ is even or $w_{q}(s)=s$.\\
2)\hspace*{1em}Either $w_{p}(s)=\overline{s}$ or $w_{pq}(s)=\overline{s}$.
\end{tabular}
%\end{center}
\\

This gives us four cases. The case $\ell (s)$ is even and 
$w_{p}(s)=\overline{s}$ is
impossible by \cite[Theorem 5.6]{JL1} since $p$ and $q$ are both odd.
Furthermore if $w_q(s) = s$ and $w_{pq}(s)=\overline{s}$,
then we also have $w_{p}(s)=\overline{s}$.  Hence we are reduced to
considering two cases. Let $\tilde{s}$ be an edge of $\Delta$ lying above
$s$.

We first consider the case $\ell (s)$ is even and $w_{pq}(s)=\overline{s}$.
Let $\gamma_{pq}$ be an element of $\cD$ representing $w_{pq}$.
Modifying $\gamma_{pq}$ by an element of $\cD^\times$ we may assume
that $\N (\gamma_{pq})=pq$ and $\gamma_{pq}\tilde{s}=\overline{\tilde{s}}$.
As in the proof of \cite[Theorem 5.6]{JL1} it follows that
$\Bpint\cong B(-1,-pq)$.  The same argument read in reverse shows
that if $\Bpint\cong B(-1,-pq)$ we can get such an edge $\tilde{s}$
which then descends to an edge $s$ of $G$ satisfying $\ell (s )=2$
and $w_{pq}(s)=\overline{s}$. The edge $s$ of $G$ in turn descends to
an edge $r$ of $\Gq$ such that $\ell (r)$ is even and $w_{p}(r)=
\overline{r}$.

Now suppose $w_{q}(s)=s$ and $w_{p}(s)=
\overline{s}$. We can find elements $\gamma_{p},\, \gamma_{q}\in
\cD$ representing $w_{p}, \,w_{q}$ respectively
with $\N (\gamma_{p})=p$ and $\N (\gamma_{q})=q$ such that
$\gamma_{p}\tilde{s}=\overline{\tilde{s}}$ and
$\gamma_{q}\tilde{s}=\tilde{s}$.
Therefore some powers of $\gamma_{p}$ and of $\gamma_{q}$ are in
$\bbQ^{\times}$. Since $\gamma_p$ and $\gamma_q$ are quadratic over
$\bbQ$, it follows $\gamma_{p}^{2}\in\bbQ^{\times}$ and
$\gamma_{q}^{2}\in\bbQ^{\times}$ (see the analysis in
\cite[Proposition 4-4]{Kur}). Hence
$\gamma_{p}^{2}=-p$ and $\gamma_{q}^{2}=-q$.  Since $\gamma_{q}$
and
$\gamma_{p}\gamma_{q}\gamma_{p}^{-1}$ both fix $\tilde{s}$, it follows that
$\gamma_{p}\gamma_{q}\gamma_{p}^{-1}=-\gamma_{q}$ and hence
$\Bpint \cong B(-p,-q)$.  Conversely suppose $\Bpint \cong B(-p,-q)$.
Pick $\gamma_{p},\,\gamma_{q}$ satisfying $\gamma_{p}^{2}=-p$,
$\gamma_{q}^{2}=-q$, and $\gamma_{p}\gamma_{q}=-\gamma_{q}\gamma_{p}$.
Let $\tilde{s}$ be the unique edge of $\Delta$  reversed by $\gamma_{p}$.
Then $\gamma_{q}(\tilde{s})$ must be $\tilde{s}$, so the image $r$
of $\tilde{s}$ in $\Gq$ is such that $\ell (r)$ is even and
$w_{p}(r)=\overline{r}$. This
completes the proof of statement (3) of Theorem \ref{local}.
\end{pf}
\section{Proof of Theorem \ref{maintheorem}}
We now prove Theorem \ref{maintheorem}.
With $B$, $p$, and $q$ satisfying the hypotheses of the theorem,
let $g(V_{B})$ be the genus of $V_{B}$. Then by the Riemann-Hurwitz
formula the genus $g$ of $\Vp$ is given by
\begin{equation}
\label{genus}
g = \frac{g(V_{B})+1}{2} - \frac{e(p)}{4}\ ,
\end{equation}
where $e(p)$ is the number of fixed points of the
$w_{p}$ on $V_{B}(\bbC )$,
cf.\ \cite[Eq.\ 3]{Ogg}. Now $e(p)$ is equal to the number
of points on $V_{B}$ with complex multiplication by $\bbZ [\sqrt{-p}]$,
which is  twice the class number $h_{p}$  of $\bbQ (\sqrt{-p} )$
if $\bbQ (\sqrt{-p})$ splits $B$ and $0$ otherwise.
As $q$ is inert in $\bbQ (\sqrt{-p})$ we  obtain:
\begin{equation}
\label{congruence}
e(p)=2h_{p}\equiv 4\!\mod 8
\end{equation}
since $p\equiv 5\! \mod 8$; see, e.g., \cite[p. 250, Exercise 25]{BS}.
In particular then $e(p)/4$ is odd.  Under our
hypotheses Eichler's mass formula gives
\begin{equation}
\label{mass}
\frac{g(V_{B})+1}{2}=1+\frac{(p-1)(q-1)-16}{24}\equiv 1 \!\mod 2
\end{equation}
(see, e.g., \cite[p. 46]{JL2}).  Combining then (\ref{genus}), 
(\ref{congruence}),
and (\ref{mass}), we see that $g$ is even.

From Theorem \ref{local} we have that $\Pic^{1}\Vp (\Qp )^{+}\neq \e$ and
$\Pic ^{1}\Vp (\bbR )^{+}\neq \e$
since $\bbQ (\sqrt{p})$ splits $B$.
Furthermore
$\Bqint\not\cong B(-1,-pq)$ since $(-1,-pq)_{p}=+1$ and
 $\Bqint\not\cong B(-p,-q)$ since $(-p,-q)_{q}= -1$.
Hence $\Pic^{1}\Vp (\Qq )^{+}=\e$ by Theorem \ref{local}(3).
Since $\Pic ^{2} V_{B}(\bbQ_{\ell} )^{+}\neq \e$ for all primes $\ell\leq \infty$
and $\Pic^{1}V_{B}(\bbQ_{\ell} )^{+}\neq \e$ for primes $\ell\neq p,q,\infty$
by \cite[Theorem 2]{JL2},
the same holds for $\Vp$.
Thus $\Pic^{g-1}(\Vp )(\Ql )^{+}=\e$ for exactly one prime $\ell\leq\infty$,
namely  $\ell=q$.  By the criterion of
Poonen and Stoll, it follows  that the jacobian of $\Vp/\bbQ$ is odd.

 It is routine that at most finitely many of the  $\VBp/\bbQ$\/'s
in Theorem \ref{maintheorem} are hyperelliptic, cf.\
\cite[Sect.\ 5]{Ogg}. The curves $\VBp$ and
$V_{B}$
have good reduction $\tVBp/\Ftwo$ and $\tVB/\Ftwo$, respectively, at $2$.
If $\VBp$ were
hyperelliptic, we would then have that
\begin{equation}
\label{1}
\#\tVp (\Ffour )\leq 2\cdot\#\Pone (\Ffour)= 10\ .
\end{equation}
But the supersingular locus on $\tVB/\Ftwo$ is rational over $\Ffour$
and its cardinality is the class number $h$ of the definite
quaternion algebra of reduced discriminant $2pq$. By Eichler's formula
(see \cite[(4-1)]{Kur}) however, we have
\begin{equation}
\label{2}
\#\tVB (\Ffour )\geq h\geq \frac{(p-1)(q-1)}{12}\ .
\end{equation}
Since $\tVB$ is a double cover of $\tVBp$ over $\Ftwo$, this  implies
that
\begin{equation}
\label{3}
\#\tVBp (\Ffour )\geq \frac{(p-1)(q-1)}{24}\ .
\end{equation}
From (\ref{1}), (\ref{2}), and (\ref{3}) we see that if
$\VBp$ is hyperelliptic, then
$(p-1)(q-1)\leq 240$.

\end{document}